\documentclass[reqno,11pt]{amsart}
\usepackage{amscd,amsfonts,amsmath,amssymb,amstext,amsthm}
\usepackage[utf8]{inputenc}
\usepackage[T2A]{fontenc}
\usepackage{cite}
\usepackage[unicode]{hyperref}
\usepackage{color}
\usepackage{xcolor}
\usepackage{comment}
\usepackage{tikz}
\usepackage{tikz-cd}
\usepackage[normalem]{ulem}

\usetikzlibrary{matrix,arrows,decorations.pathmorphing}

\makeatletter
\def\@settitle{\begin{center}%
    \baselineskip14\p@\relax
    \bfseries
    \MakeUppercase{\@title}
  \end{center}%
}
\makeatother

\newtheorem{theorem}{Theorem}
\newtheorem{lemma}{Lemma}
\newtheorem{proposition}{Proposition}
\newtheorem{corollary}{Corollary}
\theoremstyle{remark}
\newtheorem{remark}{Remark}
\newtheorem{example}{Example}
\theoremstyle{definition}
\newtheorem{definition}{Definition}
\def\C{{\mathbb C}}

\def\Cn*{{\C^n}^*}
\def\CN*{{\C^N}^*}

\def\R{{\mathbb R}}
\def\Z{{\mathbb Z}}

\def\vol{{\rm vol}}
\newcounter{par}
\newcounter{spr}
\begin{document}
\title{
How many roots of a system of random trigonometric polynomials are real?
}
\author{
B. Kazarnovskii}
\address {\noindent Institute for Information Transmission Problems, Moscow
\newline
{\it kazbori@gmail.com}.}
\begin{abstract}

The expected number of zeros of a random real polynomial of degree $N$ asymptotically equals $\frac{2}{\pi}\log N$.
On the other hand, the average fraction of real zeros of a random trigonometric polynomial of increasing degree $N$
converges to not $0$ but to $1/\sqrt 3$.
An average number of roots of a system of random trigonometric polynomials in several variables
is equal to the mixed volume of some ellipsoids depending on the degrees of polynomials.
Comparing this formula
  with Theorem BKK we prove that
the pheno\-menon of nonzero fraction of real roots remains valid.
\end{abstract}
\maketitle
%
%
%
%
%
%
\paragraph{\textbf{1. Introduction.}}
%
%
For $\lambda=(\lambda_1,\ldots,\lambda_n)\in\Z^n$ we use the notation $z^\lambda=z_1^{\lambda_1}\cdots z_n^{\lambda_n}$.
Let $\Lambda$ be a finite subset of $\Z^n$.
We remind that a function
$\sum_{\lambda\in\Lambda}a_\lambda z^\lambda$ on a complex torus $(\C\setminus0)^n$
is called a Laurent polynomial with support $\Lambda$.
The convex hull ${\rm conv}(\Lambda)$
of $\Lambda$
is called a Newton polytope of Laurent polynomial.
\begin{definition}\label{dfLorReal}
Let the set $\Lambda$ be centrally symmetric.
The Laurent polynomial $\sum_{\lambda\in\Lambda}a_\lambda z^\lambda$ is called \emph{a real Laurent polynomial},
if   $$\forall \lambda\in\Lambda\colon\: a_{-\lambda}=\overline{a_\lambda}.$$
\end{definition}
  Any real Laurent polynomial is a real function  on the compact subtorus $T^n$ of the torus
$(\C\setminus0)^n$.
\begin{definition}\label{dfTrig}
Any real Laurent polynomial with a support $\Lambda$
as a function on $T^n\subset(\C\setminus0)^n$
is called \emph{a trigonometric polynomial} on the torus $T^n$ with the support $\Lambda$.
\end{definition}
 Trigonometric polynomials with support $\Lambda$
form a real vector space $V(\Lambda)$ of dimension
equal to the number of points in $\Lambda$.
In what follows, we always consider $V(\Lambda)$
as a subspace of the space $L^2(T^n,d\chi) $,
where $d\chi$ is the normalized Haar measure on the torus $T^n$.

Let $l$ be a linear functional in the space $\Z^n\otimes_\Z\R$.
Assume that $l$ is nonzero at nonzero points of the set $\Lambda$.
We denote by $\Lambda_+$ the intersection of $\Lambda$ with
the half-space $l\geq0$.
For $0\ne\lambda\in\Lambda_+$ we put
\begin{equation}\label{eqbasis}
  \tau_\lambda(\theta)=\sqrt{2}\cos\langle\theta,\lambda\rangle,\, \tau_{-\lambda}(\theta)=\sqrt{2}\sin\langle\theta,\lambda\rangle.
\end{equation}
If $0\in\Lambda$ then we put $\tau_0(\theta)=1$.
The functions $\{\tau_\lambda\colon\lambda\in\Lambda\}$ form an orthonormal basis of the space $V(\Lambda)$.

Zeroes of trigonometric polynomials in one variable have been considered in many publications;
see \cite{ADG} and the references therein.
It turned out that, in contrast to the algebraic case
(see \cite{Ka} as well as \cite{EK} and the bibliography therein),
the expected fraction of real zeros of a random trigonometric polynomial
 of increasing degree \emph{does not tend to zero}.
Namely, the following is true.
Let $\sum_{0\leq i\leq m}a_iz^i+\overline{a_i}z^{-i}$ be a random real trigonometric polynomial.
Denote by $\mathfrak M (m)$ the average number of its zeroes located on the unit circle.
Then $\lim_{m\to\infty}\frac{1}{2m}\mathfrak M (m)=\frac{1}{\sqrt 3}$;
see Example \ref{ex12}.

It will be shown below that a similar phenomenon holds true
for trigonometric polynomials in several variables.
For this we use two theorems on the numbers of roots of systems of equations.
The first theorem states that the average number of roots of systems of random trigonometric polynomials
$\{f_1=\ldots=f_n=0\colon f_i\in V(\Lambda_i)\} $
is equal to the multiplied by $n!$ mixed volume of some ellipsoids ${\rm ell}(\Lambda_i)$,
depending on the supports $\Lambda_i$.
This theorem is a specialization of \cite[Theorem 1]{AK}
to the case $X=T^n$ and $V_i=V(\Lambda_i)$;
see also \cite{ZK}.
The second theorem states that
for almost all tuples of $n$ real Laurent polynomials with supports $\Lambda_1,\ldots,\Lambda_n$
the number of their common zeros in the complex torus $(\C\setminus0)^n$ is the same and
is equal to the mixed volume of Newton polyhedra ${\rm conv}(\Lambda_i)$ multiplied by $n!$.
This statement is a consequence of Theorem BKK; see \cite{B}.
\par\smallskip
\paragraph{\textbf{2. The average number of roots.}}
%
%
Let $\Lambda_1,\ldots,\Lambda_n$ be finite centrally symmetric sets in $\Z^n$,
and let $V(\Lambda_i)$ be a space of trigonometric polynomials with the support $\Lambda_i$,
$\dim V(\Lambda_i)=N_i=\#\Lambda_i$.
We denote by $S_i\subset V(\Lambda_i)$ the sphere of radius $1$ centered in $0$.
Let $ds_i$ be an orthogonally invariant volume form on $S_i$,
such that $\int_{S_i}ds_i=1$
(we consider the spaces $V(\Lambda_i)$ with the metric $L^2(T^n,d\chi)$).
Consider the points $s_1\in S_1,\ldots,s_n\in S_n$
as trigonometric polynomials on $T^n$,
and denote by $N(s_1,\ldots,s_n)$ the number of
isolated points of the set
of their common zeroes.
\begin{definition}\label{dfAverage}
We call
\begin{equation}\label{eqAverage}
\mathfrak M(\Lambda_1,\ldots,\Lambda_n)=\int_{S_1\times\ldots\times S_n} N(s_1,\ldots,s_n)\:ds_1\ldots ds_n.
\end{equation}
 the average number of roots of a system of random trigonometric polynomials with the supports $\Lambda_1,\ldots,\Lambda_n$.
\end{definition}
Let $\Theta_i\colon T^n\to V(\Lambda_i)$ be a map, such that
\begin{equation}\label{eqTheta}
 \forall f\in V(\Lambda_i)\colon\:\langle\Theta_i(\theta),f\rangle=\frac{1}{\sqrt N_i}f(\theta).
\end{equation}
\begin{lemma}\label{lmM1}
It is true that $\Theta_i(T^n)\subset S_i$.
\end{lemma}
\begin{proof}
Using the orthonormal basis of the space $V(\Lambda_i)$ defined in (\ref{eqbasis}),
we write the mapping $\Theta_i$ as
\begin{equation}\label{eqTheta2}
  \Theta_i(\theta)=\frac{1}{\sqrt{N_i}}\sum_{\lambda\in\Lambda}\tau_\lambda(\theta) \tau_\lambda.
\end{equation}
The statement follows from the identity $\cos^2+\sin^2=1$.
\end{proof}
Let $\mathfrak h$ and $\mathfrak h^* $ be the tangent and cotangent spaces, respectively, at the zero point of the torus $T^n$.
We denote by $\Z^n $ the integer lattice in the space $\mathfrak h$,
and by ${\Z^n}^*$ -- the dual lattice in the space $\mathfrak h^*$.
We further assume that $\Lambda_i\subset{\Z^n}^*$.
Let $F_i$ be a quadratic form in $\mathfrak h$
which is a pull-back of the
metric form on the sphere $S_i$
under the mapping  $\Theta_i$
at the zero point of $T^n$.
\begin{lemma}\label{lmM2}
It is true that $F_i(\xi)=\frac{1}{N_i}\sum_{\lambda\in\Lambda_i}\lambda_i^2(\xi)$.
\end{lemma}
\begin{proof}
Follows from (\ref{eqTheta2}).
\end{proof}
Recall that the function $h(x)=\max_{a \in A} a(x)$ in $\R^n$ is called
a support function of a compact convex body $A\subset{\R^n}^*$.
\begin{lemma}\label{lmM3}
The function
$$
h_i(\xi)=\sqrt{F_i(\xi)}=\sqrt{\frac{1}{N_i}\sum_{\lambda\in\Lambda_i}\lambda_i^2(\xi)}
$$
is a support function of some ellipsoid ${\rm ell}(\Lambda_i)$ in $\mathfrak h^*$.
If the set $\Lambda_i$ is not contained in any proper subspace of $\mathfrak h^*$,
then $\dim {\rm ell}(\Lambda_i)=n$.
\end{lemma}
\begin{proof}
For any nonnegative quadratic form $g$, it is true that
the function $\sqrt g$ is a support function of some ellipsoid.
If the form $g$ is non-degenerate, then this ellipsoid is full-dimensonal.
\end{proof}
\begin{theorem}\label{thmMixed}
It is true that
$$
\mathfrak M(\Lambda_1,\ldots,\Lambda_n)=n!\:\vol({\rm ell}(\Lambda_1),\ldots,{\rm ell}(\Lambda_n)),
$$
where in the measurement of the mixed volume it is assumed that
the volume of the fundamental cube of the lattice ${\Z^n}^*$ equals  $1$.
\end{theorem}
\begin{proof}
For $X=T^n$
(taking into account the orthogonality of the action of $T^n$ on $V(\Lambda_i)$),
the statement is easy to deduce from \cite[Теорема 1]{AK}.
\end{proof}
Note,
that the theorem \ref {thmMixed} implies the following inequalities for the mean numbers of roots.
\begin{corollary}\label{corHodge1}
For any supports $\Lambda_1,\ldots,\Lambda_n$
it is true that
%
$$
  \mathfrak M^2(\Lambda_1,\ldots,\Lambda_n)\geq\mathfrak M(\Lambda_1,\ldots,\Lambda_{n-1},\Lambda_{n-1})\cdot\mathfrak M(\Lambda_1,\ldots,\Lambda_n,\Lambda_n),
$$
$$
{\mathfrak M}^n(\Lambda_1,\ldots,\Lambda_n) \ge \mathfrak M(\Lambda_1)\cdot\ldots\cdot \mathfrak M(\Lambda_n),
$$
where, by definition, $\mathfrak M(\Lambda)=\mathfrak M(\Lambda,\ldots,\Lambda)$.
\end{corollary}
\begin{proof}
Follows from Alexandrov-Fenchel inequality (см. \cite{Al}) for mixed volumes of ellipsoids ${\rm ell}(\Lambda_i)$.
\end{proof}
\begin{remark}
The inequalities for the numbers of roots from Сorollary \ref{corHodge1} are similarly to the Hodge inequalities for the intersection indices
of projective algebraic surfaces.
The connection between the Hodge inequalities and the Alexandrov-Fenchel inequalities was independently found by A. Khovansky and B. Teissier;
see, for example, \cite{KK}.
\end{remark}
\par\smallskip
\paragraph{\textbf{3. Newton polytopes and ellipsoids.}}
%
%
Recall that with a finite centrally symmetric set $\Lambda\subset{\Z^n}^*\subset\mathfrak h^*$
we are considering two convex bodies:
the Newton polytope ${\rm conv}(\Lambda)$
and ellipsoid ${\rm ell}(\Lambda)$ with the support function
\begin{equation}\label{eq_hTrig}
  h_\Lambda(\xi)=\sqrt{\frac{1}{N}\sum_{\lambda\in\Lambda} \lambda^2(\xi)},
\end{equation}
where $N$ is the number of elements of the set $\Lambda$.
\begin{proposition}\label{pr<1}
For almost all tuples $(P_1,\ldots,P_n)$ of real Laurent polynomials with supports $\Lambda_1,\ldots,\Lambda_n$
the number of solutions of the system $P_1=\ldots=P_n=0$ is equal to $n!\:\vol({\rm conv}(\Lambda_1),\ldots,{\rm conv}(\Lambda_n))$.
\end{proposition}
\begin{proof}
The set of such tuples $(P_1,\ldots,P_n)$ is an everywhere Zariski dense subset in the space
of all tuples of Laurent polynomials with supports $\Lambda_1,\ldots,\Lambda_n$.
Therefore, the statement is a consequence of Theorem BKK.
\end{proof}
\begin{lemma}\label{lm<} $\vol({\rm ell}(\Lambda))\leq \vol({\rm conv}(\Lambda))$
\end{lemma}
\begin{proof}
From Proposition \ref{pr<1} and from Theorem \ref{thmMixed} it follows
that the right-hand side and the left-hand side of the formula
are equal to the average number, respectively, of real
and complex (i.e., all) roots of identical systems of equations.
\end{proof}
This lemma is consistent with the following geometric statement.
\begin{proposition}\label{pr<}
For any finite centrally symmetric $\Lambda$, it is true that ${\rm ell}(\Lambda)\subset{\rm conv}(\Lambda)$.
\end{proposition}
\begin{proof}
For any finite set $\Lambda$ the support function of ${\rm conv}(\Lambda)$ at the point $\xi\in\mathfrak h$
is equal to $\max_{\lambda\in\Lambda}\lambda(\xi)$.
  If $\Lambda$ is centrally symmetric
  then $\max_{\lambda\in\Lambda}\lambda(\xi)=\max_{\lambda\in\Lambda}\vert \lambda(\xi)\vert$.
  The support function of ellipsoid ${\rm ell}(\Lambda)$
  is equal to the root mean square of numbers $\{\vert \lambda(\xi)\vert\colon \lambda\in\Lambda\}$.
  The root mean square of a finite set of non-negative numbers
   does not exceed the maximum of them.
\end{proof}
\begin{theorem}\label{thmMixedTrigPart}
Let ${\rm real}(\Lambda_1,\ldots,\Lambda_n)$
be an average fraction of real roots of systems of trigonometric polynomials with
supports $\Lambda_1,\ldots,\Lambda_n$.
Then
$$
{\rm real}(\Lambda_1,\ldots,\Lambda_n)=\frac{\vol\left({\rm ell}(\Lambda_1),\ldots,{\rm ell}(\Lambda_n)\right)}{\vol\left({\rm conv}(\Lambda_1),\ldots,{\rm conv}(\Lambda_n)\right)}
$$
\end{theorem}
\begin{proof}
Follows from Theorem \ref{thmMixed} and proposition \ref{pr<1}.
\end{proof}
\begin{example}\label{ex11}
In the case $n=1$, the segments ${\rm ell}(\Lambda)$ and ${\rm conv}(\Lambda)$
coincide, if and if $\Lambda=\{\lambda,-\lambda\}$.
From the Theorem \ref{thmMixedTrigPart} it follows 
that all zeros of any polynomial of the form $az^\lambda+\bar az^{-\lambda}$
belong to the unit circle,
which is true, because these zeros are the $(2\lambda)$-roots of $\bar a/a$.
\end{example}
\begin{example}\label{ex12}
If $\Lambda=\{-\lambda,\ldots,-1,0,1,\ldots,\lambda\}$,
then the ellipsoid ${\rm ell (\Lambda)}$
is a segment
$$\left[-\sqrt{\frac{\lambda(\lambda+1)}{3}},\sqrt{\frac{\lambda(\lambda+1)}{3}}\right]$$
Hence
Theorem \ref{thmMixedTrigPart} implies that,
for large $\lambda$,
the average fraction of real zeros of trigonometric polynomials with support $\Lambda$
tends to $\sqrt{\frac {1}{3}}\simeq 0.577$.
\end{example}
\par\smallskip
\paragraph{\textbf{4. Asymptotics of the average fraction of real roots.}}
Let $\Delta\subset\mathfrak b^*$ be a centrally symmetric convex body, and $\Lambda=\Delta\cap{\Z^n}^*$.
Further, we assume
that the following condition (*) is satisfied

\textbf{(*)}\ \
\emph{if $\dim\Delta=k<n$, then $\Delta$
is contained in some $k$-dimensional subspace of $\mathfrak h^*$,
generated by the vectors of the lattice ${\Z^n}^*$.}

For an integer $m>0$ we put
$\Delta_m=m\Delta$, $\Lambda_m=\Delta_m\cap{\Z^n}^*$, $N_{\Lambda,m}=\#\Lambda_m$.
Recall
that
a function $ h_\Lambda=\sqrt{F_\Lambda}$,
  where
$$
F_\Lambda(\xi)=\frac{1}{N_\Lambda}\sum_{\lambda\in\Lambda}\lambda^2(\xi),
$$
is the support function of ellipsoid
${\rm ell}(\Lambda)\subset\mathfrak h^*$.
\begin{lemma}\label{lmLim}
For $m\to\infty$, the sequence of functions $\frac{1}{m^2}F_{\Lambda,m}$
converges locally uniformly to some quadratic form $F_\Delta$.
\end{lemma}
\begin{proof}
Choose a basis of the lattice ${\Z^n}^*$ and declare it
an orthonormal basis of the space $\mathfrak h^*$.
Using the appropriate metric,
we identify the spaces $\mathfrak h$ and $\mathfrak h^*$,
as well as integer lattices $\Z^n\subset\mathfrak h$ and ${\Z^n}^*\subset\mathfrak h^*$.
Let $\dim\Delta=k$.
Then
$$
  \frac{1}{m^2} F_{\Lambda,m}(\xi)=\frac{1}{m^2N_{\Lambda,m}}\sum_{\lambda\in\Lambda_m}\langle\lambda,\xi\rangle^2=
  \frac{m^k}{N_{\Lambda,m}}
  \sum_{\alpha\in\frac{\Lambda_m}{m}}
  \langle\alpha,\xi\rangle^2\frac{1}{m^k}.
$$
Now notice
that $N_{\Lambda,m}\asymp\vol_k(\Delta)m^k$,
and $\sum_{\alpha\in\frac{\Lambda_m}{m}}
  \langle\alpha,\xi\rangle^2\frac{1}{m^k}$
   is the integral sum
for an integral of the function $f(x)=\langle x,\xi\rangle^2$ over the polytope $\Delta$.
Hence
\begin{equation}\label{eqEllSupp}
  \frac{1}{m^2}F_{\Lambda,m}(\xi)\to F_\Delta(\xi)=\frac{1}{\vol_k(\Delta)}\int_\Delta\langle x,\xi\rangle^2\:dx
\end{equation}
\end{proof}
\begin{corollary}\label{corLim}
Let $m\to\infty$. Then

{\rm(1)}\ the sequence of functions $\frac{1}{m}h_{\Lambda,m}$
converges locally uniformly to the support function $h_\Delta$
of some ellipsoid ${\rm ell}(\Delta)$

{\rm(2)}\
the sequence of ellipsoids $\frac{1}{m}{\rm ell}(\Lambda_m)$
converges in the Hausdorff topology to the ellipsoid ${\rm ell}(\Delta)$.
\end{corollary}
\begin{proof}
Since, by definition, $h_\Delta=\sqrt{F_\Delta}$,
then both statements are direct consequences of Lemma \ref{lmLim}.
\end{proof}
\begin{lemma}\label{lmLimConv}
If the condition {\rm(*)} is satisfied for the convex body $\Delta$,
then for $m\to\infty$ the sequence of convex polytopes $\frac{1}{m}{\rm conv}(\Lambda_m)$
converges in Hausdorff topology to the convex body $\Delta$.
\end{lemma}
\begin{proof}
This follows from the definition of the set $\Lambda_m$.
\end{proof}
\begin{theorem}\label{thmMixedTrigAsympPart}
Let $\Delta_1,\ldots,\Delta_n$ be convex bodies in the space $\mathfrak h^*$,
satisfying the condition {\rm(*)},
and $\Lambda_i=\Delta_i\cap{\Z^n}^*$.
Then
\begin{equation}\label{eqMain}
  \lim_{\inf(m_1,\ldots,m_n)\to\infty}{\rm real}_n((\Lambda_1)_{m_1},\ldots,(\Lambda_n)_{m_n})=
\frac{\vol({\rm ell}(\Delta_1),\ldots,{\rm ell}(\Delta_n))}{\vol(\Delta_1,\ldots,\Delta_n)}
\end{equation}
\end{theorem}
\begin{proof}
Let $(m_1)_i,\ldots,(m_n)_i$ be an increasing sequence
of tuples of positive integers,
and let $(\Lambda_1)_{(m_1)_i},\ldots,(\Lambda_n)_{(m_n)_i}$
be the corresponding sequence of tuples of supports.
Then the passage to the limit based on the corollary \ref{corLim} (2) and the lemma \ref{lmLimConv},
gives the required statement.
\end{proof}
\begin{corollary}\label{corHom}
Let $\alpha_1,\ldots,\alpha_n$ be positive numbers.
Then, when replacing $\Delta_i\to\alpha_i\Delta_i$
(with such replacements the condition {\rm(*)} remains satisfied)
the asymptotic behavior of ${\rm real}_n(\Lambda_{m_1},\ldots,\Lambda_{m_n})$ is preserved.
\end{corollary}
\begin{proof}
From (\ref{eqEllSupp}) it follows, что $h_{\alpha\Delta_i}=\alpha h_{\Delta_i}$.
Therefore, the numerator and denominator of a fraction in right-hand side (\ref{eqMain}) 
are multiplied by $\alpha_1\cdots\alpha_n$.
\end{proof}
\paragraph{\textbf{5. An example of calculating the asymptotics of the average fraction of real roots.}}
%
%
Here we calculate the limit
$$
\lim_{\inf(m_1,\ldots,m_n)\to\infty}{\rm real}_n(\Lambda_{m_1},\ldots,\Lambda_{m_n})
$$
for balls $\Delta_i$ radii $R_i$.
\begin{theorem}\label{thmLimPart}
It is true that
$$
\lim_{\inf(m_1,\ldots,m_n)\to\infty} {\rm real}_n(m_1,\ldots,m_n)=\left(\frac{\beta_n}{\sigma_n}\right)^{\frac{n}{2}},
$$
where $\beta_n=\int_{-1}^1x^2(1-x^2)^{\frac{n-1}{2}}dx$, and $\sigma_n$ is a volume of $n$-dimensional ball of radius $1$.
\end{theorem}
We give the values of $\beta_n$ for $n\leq20$
(if $n=1$, then $\sqrt{\frac{1}{2}\beta_1}=\frac{1}{\sqrt 3}$; see Example \ref{ex12}).
\noindent
\begin{table}[h!]
  \begin{center}
   \begin{tabular}{r| l l l l l l l l l l}
       {\scriptsize $n$}&{\scriptsize $1$}&{\scriptsize $2$}&{\scriptsize $3$}&{\scriptsize $4$}&{\scriptsize $5$}&{\scriptsize $6$}&{\scriptsize $7$}&{\scriptsize $8$}&
       {\scriptsize $9$}&{\scriptsize $10$}\\
{\scriptsize $\beta_n$}&
      {\scriptsize $\frac{2}{3}$}&
      {\scriptsize $\frac{\pi}{8}$}&
      {\scriptsize $\frac{4}{15}$}&
      {\scriptsize  $\frac{\pi}{16}$}&
      {\scriptsize $\frac{16}{105}$}&
      {\scriptsize $\frac{5\pi}{128}$}&
      {\scriptsize $\frac{32}{315}$}&
      {\scriptsize $\frac{7\pi}{256}$}&
      {\scriptsize $\frac{256}{3465}$}&
      {\scriptsize $\frac{21\pi}{1024}$}
      \\
      \hline
       {\scriptsize $n$}
       &{\scriptsize $11$}&
       {\scriptsize $12$}&
       {\scriptsize $13$}&
       {\scriptsize $14$}&
       {\scriptsize $15$}&
       {\scriptsize $16$}&
       {\scriptsize $17$}&
       {\scriptsize $18$}&
       {\scriptsize $19$}&
       {\scriptsize $20$}
                     \\
{\scriptsize $\beta_n$}&
      {\scriptsize $\frac{512}{9009}$}&
      {\scriptsize $\frac{33\pi}{2048}$}&
      {\scriptsize $\frac{4096}{109395}$}&
{\scriptsize $\frac{429\pi}{32768}$}&
{\scriptsize $\frac{2048}{45045}$}&
{\scriptsize $\frac{715\pi}{65536}$}&
{\scriptsize $\frac{65536}{2078505}$}&
{\scriptsize  $\frac{2431\pi}{262144}$} &
{\scriptsize $\frac{131072}{4849845}$}&
{\scriptsize $\frac{4199\pi}{524288}$}
     \\
     \end{tabular}
  \end{center}
\end{table}
\begin{remark}
The expression $x^2(1-x^2)^{\frac{n-1}{2}} dx$ is the so-called \emph{differential Tche\-bichef binomial}.
Tchebichef proved \cite{Ch} that the binomial $x^m(a+bx^n)^p dx$ is non-integrable outside
three (known to L. Euler) cases.
For odd $n$, the above expression refers to the first one,
and for even $n$ - to the third case.
\end{remark}
Application of Lemma \ref{lmLimConv} and Theorem \ref{thmMixedTrigAsympPart}
reduces the proof of Theorem \ref{thmLimPart} to the following statement.
\begin{proposition}\label{prLimBeta} It is true that
$
{\rm ell}(\Delta_i)=\sqrt{\frac{\beta_n}{\sigma_n}}\:\Delta_i.
$
\end{proposition}
\begin{proof}
According to Corollary \ref{corLim},
the sequence of ellipsoids $\frac{1}{m}{\rm ell}\left((\Lambda_i)_m\right)$
converges to the ellipsoid ${\rm ell}(\Delta_i)$
with support function
$$
h_{\Delta_i}(\xi)=\sqrt{\frac{1}{\vol_k(\Delta_i)}\int_{\Delta_i}\langle x,\xi\rangle^2\:dx}.
$$
The right side of the equation,
as it is easy to see,
is equal to  $|\xi|R_i\sqrt{\frac{\beta_n}{\sigma_n}}$,
where $R_i$ is the radius of the ball $\Delta_i$.
Thus,
$h_{\Delta_i}(\xi)$ is the support function of the ball of radius $R_i\sqrt{\frac{\beta_n}{\sigma_n}}$
centered at the point $0$.
This ball is $\sqrt{\frac{\beta_n}{\sigma_n}}\Delta_i$.
The proposition is proved.
\end{proof}
\end{document}